\documentclass{amsart}

\usepackage{amssymb}
\usepackage{graphicx}
\usepackage{flafter}

\newtheorem{theorem}{Theorem}[section]

\newtheorem{lemma}[theorem]{Lemma}
\newtheorem{cor}[theorem]{Corollary}

\theoremstyle{definition}

\theoremstyle{remark}

\numberwithin{equation}{section}

\newcommand\blfootnote[1]{%
  \begingroup
  \renewcommand\thefootnote{}\footnote{#1}%
  \addtocounter{footnote}{-1}%
  \endgroup
}

\begin{document}

\title[Exponentially integrable distortion and pointwise rotation]{Pointwise rotation for mappings with exponentially integrable distortion}
\keywords{Mappings of finite distortion, rotation, exponentially integrable distortion. \newline \indent The author was financially supported by the Jenny and Antti Wihuri Foundation and by The Centre of Excellence in Analysis and Dynamics Research (Academy of Finland, decision 271983)}

\author{Lauri Hitruhin}
\address{University of Helsinki, Department of mathematics and Statistics, P-O. Box 68, FIN-00014 University of Helsinki, Finland}
\email{lauri.hitruhin@helsinki.fi}


\begin{abstract}
We prove an upper bound for pointwise rotation of mappings with $p$-exponentially integrable distortion. We also show that this bound is essentially optimal by providing examples which attain this rotation up to a constant multiplication. 
\end{abstract}

 \maketitle

\section{Introduction}\label{INTRO}
Let $f: \mathbb{C}\to \mathbb{C}$ be a $K$-quasiconformal mapping normalized by $f(0)=0$ and $f(1)=1$. Then the pointwise stretching properties of this mapping are captured by the classical H\"older continuity result 
\begin{equation}\label{venytyskvaseille}
\frac{1}{c_K}|z|^{K}\leq |f(z)|\leq c_K|z|^{\frac{1}{K}},
\end{equation}
for all $0\leq |z|\leq 1$. On the other hand, the rotational properties of quasiconformal mappings, and also of mappings of finite distortion, have been earlier studied (see, for example, \cite{BFP}, \cite{GM} and \cite{J}) by restricting to mappings between annuli and then measuring the maximal rotation of the inner circle. Recently in \cite{AIPS} Astala, Iwaniec, Prause and Saksman proposed a new pointwise approach to rotation of quasiconformal mappings, dropping assumptions regarding annuli, and proved that
\begin{equation}\label{kiertokvaseille}
|\arg(f(r))|\leq \frac{1}{2}\left( K-\frac{1}{K} \right)\log\left( \frac{1}{r} \right)+c_K,
\end{equation}
where we use the same normalization as in \eqref{venytyskvaseille}, $\arg$ is the principal branch of the argument and $0<r<1$. \newline \newline
In this paper we study mappings of finite distortion with $p$-exponentially integrable distortion function, that is 
\begin{equation*}
e^{pK_f(z)}\in L^{1}_{\text{loc}}(\mathbb{C}).
\end{equation*}
For this class of mappings the analogue to \eqref{venytyskvaseille} is given by the modulus of continuity results in \cite{HK} and \cite{OZ}, which show that
\begin{equation*}
e^{-\frac{c}{p}\log^{2}\left( \frac{1}{|z|} \right)} \lesssim |f(z)|  \lesssim \frac{1}{\log^{\frac{p}{2}}\left( \frac{1}{|z|} \right)},
\end{equation*}
where $0<|z|<\epsilon_{f,p}$ and $f(0)=0$. So, the pointwise stretching is well understood for mappings with exponentially integrable distortion. However, the question regarding the analogous result for the pointwise rotation \eqref{kiertokvaseille} in this class of mappings has remained open. Our aim is to answer this question in the form of the following theorem: 
\begin{theorem}  \label{Theorem1}
Fix an arbitrary $p>0$ and let $f: \mathbb{C} \to \mathbb{C} $ be a mapping of finite distortion such that $e^{pK_f(z)}\in L^{1}_{\text{loc}}(\mathbb{C})$, normalized by the conditions $f(0)=0$ and $f(1)=1$. Then 
\begin{equation} \label{Eq1}
|\arg(f(z)) |\leq \frac{c}{p}\log^{2}\left( \frac{1}{|z|} \right),
\end{equation}
when $|z| \to 0$. More precisely, for every such $f$ there exists a constant $\epsilon_f >0$ such that \eqref{Eq1} holds for all $z$ which satisfy $0<|z|<\epsilon_f$. Here $c$ is a fixed constant that does not depend on the parameter $p$ or the mapping $f$ and $\arg$ is the principal branch of the argument. 
\end{theorem}
We will also show that Theorem \ref{Theorem1} is optimal, up to the exact value of the constant $c$. We do this by providing for any given $\epsilon>0$ a mapping $h$, which satisfies the assumptions of Theorem \ref{Theorem1}, such that
\begin{equation}\label{maksimikierto}
|\arg(h(r)) | = \frac{1-\epsilon}{2p}\log^{2}\left( \frac{1}{r} \right),
\end{equation}
for every $0<r<\frac{1}{2}$. It remains open if the constant $\frac{1}{2p}$ is optimal in the inequality \eqref{Eq1}, but we are inclined to believe so. \newline \newline
Moreover, in \cite{AIPS} it was proved that given a $K$-quasiconformal mapping $f: \mathbb{C} \to \mathbb{C}$, again with the normalization $f(0)=0$ and $f(1)=1$, the pointwise stretching bounds the principal branch of the argument by
\begin{equation*}
|\arg(f(r))|\leq c_K|\log|f(r)||,
\end{equation*}
when $r\in (0,\epsilon_f)$. Using this relation between stretching and argument they defined the pointwise rotation of a mapping $f$ at a point $z_0\in \mathbb{C}$ as the limit 
\begin{equation*}
\gamma_f(z_0)=  \lim_{n \to \infty} \frac{\arg(f(z_0+r_n)-f(z_0))}{\log | f(z_0+r_n)-f(z_0) |},
\end{equation*} 
where $r_n \to 0$ is some decreasing sequence of positive radii. \newline \newline 
However, for mappings with exponentially integrable distortion the argument is in general not bounded by $c| \log|f(r)||$. This follows as the mapping
\begin{equation}\label{esimkiertoeiveny}
h_0(z)=  \left\{
  \begin{array}{l l}
    ze^{-ic_2\log^{\frac{3}{2}}\left( \frac{1}{|z|} \right)} & \quad \text{if $z\in \mathbb{D}$}\\
    
    z & \quad \text{if $z\notin \mathbb{D}$ }
  \end{array} \right. 
\end{equation}
has $p$-exponentially integrable distortion with a suitable choice for the parameter $c_2$, but clearly satisfies
\begin{equation*}
|\arg(h_0(r))|=\left| c_2\log^{\frac{3}{2}}\left( \frac{1}{r} \right)\right|> c |\log(r)|=c |\log|h_0(r)||
\end{equation*}
for any constant $c$, when the radius $r$ is small enough. Therefore we use the formulation of Theorem \ref{Theorem1}, and later when defining the pointwise rotation, use the formulation \eqref{rot}, instead of the one in \cite{AIPS}. Nevertheless, we will later on see how the pointwise stretching will bound from above the pointwise rotation even in the class of mappings with exponentially integrable distortion. \newline \newline
The mappings $h$ used to prove the optimality, up to the constant $c$, for Theorem \ref{Theorem1} are defined by
\begin{equation}\label{Esim}
h(z)=  \left\{
  \begin{array}{l l}
    \frac{z}{|z|}|z|^{c_1\log\left( \frac{1}{|z|} \right)-ic_2\log\left( \frac{1}{|z|} \right)} & \quad \text{if $z\in B\left(0,\frac{1}{2}\right)$}\\
    \frac{z}{|z|}|z|^{c_1\log(2)-ic_2 \log(2)} & \quad \text{if } z\in \mathbb{D}\setminus B\left(0,\frac{1}{2}\right) \\
    z & \quad \text{if $z\notin \mathbb{D}$ }
  \end{array} \right. 
\end{equation}
We will show that these mappings have $p$-exponentially integrable distortion if  the parameters $c_1>0$, $c_2 \in \mathbb{R}$ satisfy
\begin{equation}\label{vakioyht}
c_1+\frac{c_{2}^{2}}{c_1}<\frac{1}{p},
\end{equation}
and that for any $\epsilon>0$ we can choose the parameters $c_1$ and $c_2$ such that \eqref{maksimikierto} holds for every $0<r<\frac{1}{2}$.  \newline\newline
\section{Definitions and prerequisites }\label{PRE}
Let $\Omega \subset \mathbb{C}$ be a domain and $f: \Omega \to \mathbb{C}$ a sense preserving homeomorphism. We say that $f$ has finite distortion if the following conditions hold:
\begin{itemize}
\item $f\in W_{\text{loc}}^{1,1}(\Omega)$
\item $J_f(z)\in L^{1}_{\text{loc}}(\Omega)$
\item $|Df(z)|^2\leq J_f(z)K(z) \qquad \text{almost everywhere in $\Omega$},$
\end{itemize}
for a measurable function $K(z)\geq 1$, which is finite almost everywhere. The smallest such function is denoted by $K_f(z)$ and is called the distortion function of $f$. Here we follow the traditional notation where $Df(z)$ denotes the differential matrix of $f$ at the point $z$, $J_f(z)$ denotes the Jacobian at the point $z$ and the norm $|Df(z)|$ is defined by 
\begin{equation*}
|Df(z)|=\max \{|Df(z) e|: e\in \mathbb{C}, |e|=1\}.
\end{equation*}
Such a mapping $f$ is said to have a $p$-exponentially integrable distortion if its distortion function $K_f(z)$ satisfies
\begin{equation*}
e^{pK_f(z)} \in L_{\text{loc}}^{1}(\Omega).
\end{equation*}
For a detailed study of these mappings see, for example, \cite{AIM}. \newline \newline 
Let $f: \mathbb{C} \to \mathbb{C}$ be a mapping of finite distortion. When we study the pointwise rotation of the mapping $f$ at a point $z_0 \in \mathbb{C}$ we examine the change of the argument of $f\left( z_0+te^{i\theta} \right)-f(z_0)$ as the parameter $t$ goes from $1$ to $r>0$, which can be written as
\begin{equation*}
\left|\arg (f(z_0+re^{i\theta})-f(z_0))-\arg( f\left(z_0+e^{i\theta}\right)-f(z_0))\right|.
\end{equation*}
This can also be understood as the maximal winding of the path $f\left(\left[ z_0+e^{i\theta}, z_0+re^{i\theta} \right]\right)$ around the point $f(z_0)$. As we are interested in the maximal change of the argument, from an arbitrary direction $\theta$, we will study the supremum 
\begin{equation}\label{llll}
\sup_{\theta\in [0,2\pi)}\left|\arg (f(z_0+re^{i\theta})-f(z_0))-\arg( f\left(z_0+e^{i\theta}\right)-f(z_0))\right|.
\end{equation}
And finally, as we are interested in how fast \eqref{llll} grows in the limit case $r\to 0$, we define the pointwise rotation at the point $z_0$ by 
\begin{equation}\label{rot}
\gamma_{f}(z_0)= \limsup_{r\to 0} \frac{\sup_{\theta\in [0,2\pi)}\left|\arg (f(z_0+re^{i\theta})-f(z_0))\right|}{\log^{2}\left( \frac{1}{r} \right)},
\end{equation}
where $\arg$ is any branch of the argument and $\log^{2}\left( \frac{1}{r} \right)$ is a normalization factor, which is of the right order due to Theorem \ref{Theorem1}. Moreover, we can use interchangeably with $\arg(f(z))$ the notion $\Im \log(f(z))$, where we choose the corresponding branch of the logarithm. This is useful in some situations, for example, when calculating rotation of the mappings \eqref{Esim}.   \newline \newline 
Next we note, that for mappings with exponentially integrable distortion we can normalize general pointwise rotation in terms of Theorem \ref{Theorem1}.
\begin{cor}\label{yleinentapa}
Let $f:\mathbb{C} \to \mathbb{C}$ be a mapping of finite distortion such that $e^{pK_f(z)}\in  L^{1}_{\text{loc}}(\mathbb{C})$ and let $z_0\in \mathbb{C}$ be arbitrary. Then there exists a normalized mapping $f_0$, which satisfies the conditions of Theorem \ref{Theorem1} with the same parameter $p$, such that the pointwise rotation of the mapping $f$ around the point $z_0$ is the same as the pointwise rotation of the mapping $f_0$ around the origin. 
\end{cor}
\textit{Proof.} Define $f_0(z)=h[f(z_0+z)-f(z_0)]$,  where $h$ is the constant for which $f_0(1)=1$. It is easy to see that $f_0$ satisfies the assumptions of Theorem \ref{Theorem1} with the desired parameter $p$. Moreover, the pointwise rotation of $f_0$ around the origin is the same as the pointwise rotation of $f$ around the point $z_0$, since the constant $h$ plays no role in \eqref{rot}. \newline \newline
Hence, when studying the pointwise rotation of mappings with exponentially integrable distortion we can restrict ourselves to mappings that satisfy the assumptions of Theorem \ref{Theorem1} and measure the rotation at the origin. Then note, that given an arbitrary function $f$ satisfying the assumptions of Theorem \ref{Theorem1} the inequality \eqref{Eq1} controls the rotation \eqref{rot} and implies that 
\begin{equation*}
\gamma_f(0) \leq \frac{c}{p}.
\end{equation*}
Thus Theorem \ref{Theorem1} and Corollary \ref{yleinentapa} together prove that the limit \eqref{rot} exists and is finite for an arbitrary mapping $f:\mathbb{C}\to \mathbb{C}$ with exponentially integrable distortion and for an arbitrary point $z_0\in \mathbb{C}$. Moreover, the examples \eqref{Esim} show that the rotation 
\begin{equation*}
\frac{1-\epsilon}{2p}
\end{equation*}
is attainable for any $\epsilon>0$. \newline \newline
In the proof of Theorem \ref{Theorem1} the modulus of path families will play an important role. We give here the main definitions, but for a closer look see, for example, \cite{V}. We say that an image of a continuous mapping $\gamma: I \to \mathbb{C}$, where $I$ is an interval, is a path. We will denote both the mapping and its image by $\gamma$. Let $\Gamma$ be a family of paths. We say that a measurable function $\rho: \mathbb{C} \to [0,\infty]$ is admissible with respect to $\Gamma$ if 
\begin{equation}\label{ab}
\int_{\gamma} \rho(z) dz \geq 1,
\end{equation} 
for every locally rectifiable path $\gamma \in \Gamma$. We denote the modulus of a path family $\Gamma$ by $M(\Gamma)$ and define it by
\begin{equation}\label{modmäär}
M(\Gamma)= \inf_{\rho \enskip \text{admissible}} \int_{\mathbb{C}} \rho^2(z) dz.
\end{equation} 
We will also need a weighted version of \eqref{modmäär}, where the weight $\omega(z): \mathbb{C} \to [0,\infty]$ is measurable and locally integrable, defined by
\begin{equation*}
M_{\omega}(\Gamma)=\inf_{\rho \enskip \text{admissible}} \int_{\mathbb{C}} \rho^2(z) \omega(z) dz.
\end{equation*}
Note that in \eqref{ab} only locally rectifiable paths are considered. \newline \newline
We will also need the subsequent modulus of continuity result which follows from [\cite{HK}, Theorem B] by Herron and Koskela. Their result states, that for mappings $f$ satisfying the assumptions of Theorem \ref{Theorem1} 
\begin{equation}\label{modcont2}
|f(z)| \geq e^{-\frac{c}{p}\log^{2}\left( \frac{1}{|z|} \right)},
\end{equation}
where the constant $c$ is fixed and $0<|z|<\epsilon$, for some positive $\epsilon$ which depends on $f$ and $p$. The exact value of the constant $c$ in \eqref{modcont2} is not known, and due to this we are not able to calculate the explicit value of the constant $c$ in Theorem \ref{Theorem1}. \newline \newline
We will use $c$ to denote a generic constant which does not depend on any parameters and which value can change even on the same line of inequalities. We denote the unit disc by $\mathbb{D}$, the boundary of a disc $B(x,r)$ is denoted by $\partial B(x,r)$, the radius of a disc $B$ is denoted by $r(B)$ and $cB(x,r)=B(x,cr)$ for every $c>0$. \newline \newline
\section{Proofs}\label{Proofs}
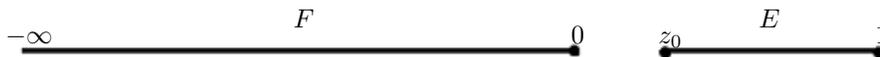
\begin{figure}
\begin{center}
\includegraphics[scale=1]{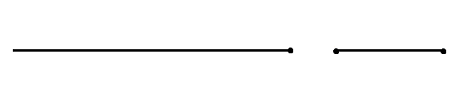}
\put(-54,40){$E$}
\put(-230,40){$F$}
\put(-92,34){$z_0$}
\put(-10,34){$1$}
\put(-125,34){$0$}
\put(-339,34){$-\infty$}
\caption{The sets $E$ and $F$.}
\label{FIG2}
\end{center}
\end{figure}

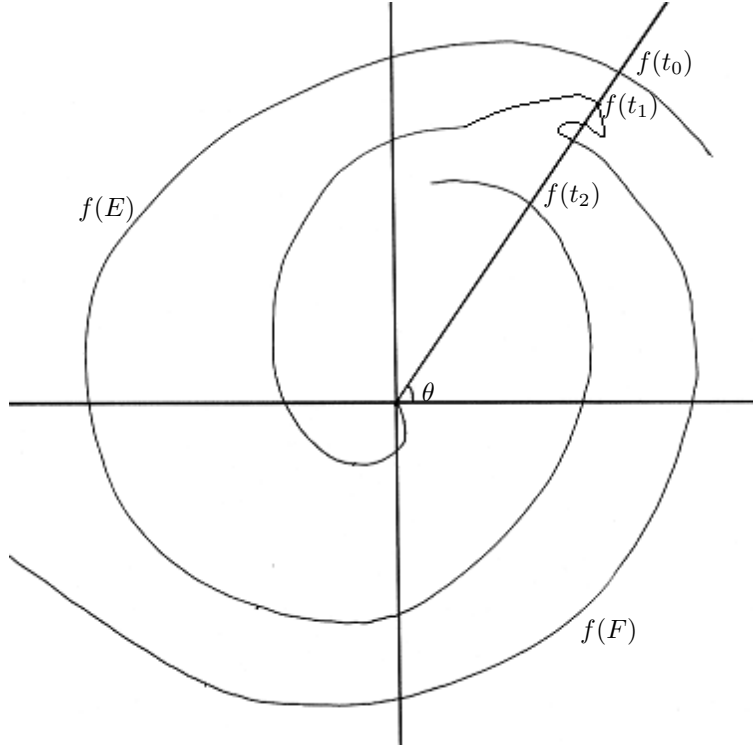
\begin{figure}
\begin{center}
\includegraphics[scale=1]{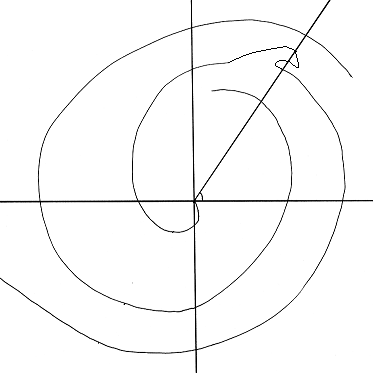}
\put(-65,40){$f(F)$}
\put(-255,200){$f(E)$}
\put(-125,130){$\theta$}
\put(-44,255){$f(t_0)$}
\put(-57,239){$f(t_1)$}
\put(-79,205){$f(t_2)$}
\caption{The paths $f(E)$ and $f(F)$ must cycle alternately around the origin.}
\label{FIG1}
\end{center}
\end{figure}

We will formulate Theorem \ref{Theorem1} with a slightly different, but clearly equivalent, normalization.
\begin{theorem}\label{Theorem2} 
Fix an arbitrary $p>0$ and let $f: \mathbb{C} \to \mathbb{C} $ be a mapping of finite distortion such that $e^{pK_f(z)}\in L^{1}_{\text{loc}}(\mathbb{C})$. Normalize it by the conditions $f(0)=0$, $f(\mathbb{D})\subset \mathbb{D}$ and $\overline{f(\mathbb{D})}\cap \partial \mathbb{D}\neq \emptyset$, and fix any branch of the argument. Then 
\begin{equation}\label{taasviite!}
|\arg(f(z)) |\leq \frac{c}{p}\log^{2}\left( \frac{1}{|z|} \right),
\end{equation}
when $|z|\to 0$. More precisely, for every such $f$ there exists a constant $\epsilon_f >0$ such that \eqref{taasviite!} holds for all $z$ which satisfy $0<|z|<\epsilon_f$. Here $c$ is some fixed constant that does not depend on the parameter $p$ or the mapping $f$. 
\end{theorem}
We will prove Theorem \ref{Theorem1} in the form of Theorem \ref{Theorem2} as the condition $|f(z)|<1$ when $z\in \mathbb{D}$ simplifies notation in the proof. \newline \newline 
Fix an arbitrary $p>0$, let $f$ be a function satisfying the assumptions of Theorem \ref{Theorem2} and let $z\in \mathbb{D}\setminus \{0\}$ be arbitrary. We will estimate $\left|\arg(f(z))-\arg\left(f\left(\frac{z}{|z|}\right)\right)\right|$, which is the winding of the path $f\left(\left[z, \frac{z}{|z|} \right]\right)$ around the origin, using the modulus of a path family. Our aim is to show that 
\begin{equation}\label{uusiviitte}
\left|\arg(f(z))-\arg\left(f\left(\frac{z}{|z|}\right)\right)\right|\leq \frac{c}{p}\log^{2}\left( \frac{1}{|z|} \right),
\end{equation}
when $0<|z|<\epsilon_f$, which is enough to show Theorem \ref{Theorem2} as $\left|\arg\left(f\left(\frac{z}{|z|}\right)\right)\right|$ is finite for any branch of the argument. \newline \newline 
To this end, we will use Corollary 4.2 from \cite{KO}, due to Koskela and Onninen, which treats capacity and moduli inequalities for a very general class of mappings satisfying certain Orlicz-type conditions. As we are interested in mappings satisfying the assumptions of Theorem \ref{Theorem2}, notably in homeomorphisms with exponentially integrable distortion, it is easy to see that the Orlicz-type conditions in their result are satisfied. Moreover, as path families $\Gamma$ we will use all paths connecting two closed separate sets, which we will specify later, and would like to get a weighted moduli inequality between $\Gamma$ and $f(\Gamma)$. These choices for path families satisfy the assumptions regarding paths made in [\cite{KO}, corollary 4.2]. Thus, applying their result to mappings that satisfy the assumptions of Theorem \ref{Theorem2}, we obtain the weighted moduli inequality  
\begin{equation}\label{A}
M(f(\Gamma))\leq M_{K_f}(\Gamma),
\end{equation}
which we shall use to prove the estimate \eqref{uusiviitte}. \newline \newline
Fix an arbitrary point $z_0\in \mathbb{D}\setminus\{0\}$. Using rotation we can assume that $0<z_0<1$, which will slightly simplify notation. Then define the line segments $E=\left[z_0,1\right]$ and $F=(-\infty,0]$, see Figure \ref{FIG2}, and let $\Gamma$ be the family of all paths joining the sets $E$ and $F$. Our intent is to estimate the values of $M_{K_f(z)}(\Gamma)$ and $M(f(\Gamma))$, and use the inequality \eqref{A} to obtain the desired upper bound for the winding of $f(E)$, which is the same as the winding in \eqref{uusiviitte}. Without loss of generality we can additionally assume that $\arg(f(z_0))-\arg(f(1))\geq 0$. \newline \newline 
\noindent Let us first estimate $M_{K_f}(\Gamma)$ from above when $z_0$ is small. To this end, construct balls $B_j=B(2^{j-1}z_0,2^{j-1}z_0)$, where $j$ goes trough numbers $1,2,...,n$ and $n$ is the smallest number for which $2^{n-1}z_0 \geq 1$. Then define 
\begin{equation*}
\rho_{0}(z)=  \left\{
  \begin{array}{l l}
    \frac{2}{z_0} & \quad \text{if $z\in B_1$}\\
    \frac{2}{2z_0} & \quad \text{if $z\in B_2 \setminus B_1$ }\\
    \vdots & \quad \vdots \\
    \frac{2}{2^{n-1}z_0} & \quad \text{if $z\in B_n \setminus B_{n-1}$ }\\
    0 & \quad \text{otherwise }
  \end{array} \right.
\end{equation*} 
To see that $\rho_{0}$ is admissible for the path family $\Gamma$ note that every point $z\in E$ belongs to some ball $\frac{1}{2}B_j$, we have $\rho_{0}(z)\geq \frac{2}{r(B_j)}$ when $z\in B_j$ and $B_j \cap F= \emptyset$ for every $j$. Using the mapping $\rho_{0}$ we estimate
 
\begin{equation}\label{turha2}
M_{K_f}(\Gamma)\leq  \int_{\mathbb{C}} K_f(z) \rho_{0}^{2}(z) dz.
\end{equation}
 To estimate this integral further we use the elementary inequality  
\begin{equation*}
ab\leq a\log(a+1)+e^{b}-1,
\end{equation*}
which holds for any $a,b\geq 0$, to obtain the pointwise inequality
\begin{equation}\label{turha1}
\frac{1}{p}(pK_f(z))\rho_{0}^{2}(z) \leq \frac{1}{p}(e^{pK_f(z)}-1)+\frac{1}{p}\rho_{0}^2(z)\log(1+\rho_{0}^2(z))
\end{equation}
that holds almost everywhere. By combining \eqref{turha2} and \eqref{turha1} we obtain 
\begin{equation}\label{olkoon1}
\begin{split}
\int_{\mathbb{C}} K_f(z)\rho_{0}^{2}(z) dz  \leq & \frac{1}{p}\int_{\mathbb{C}\setminus \{z\in \mathbb{C}: \rho_{0}(z)=0\}} e^{pK_f(z)}-1 dz  \\
& + \frac{1}{p}\int_{\mathbb{C}} \rho_{0}^{2}(z) \log(1+\rho_{0}^{2}(z)) dz .
\end{split}
\end{equation}
As $f$ is a mapping with $p$-exponentially integrable distortion and $\sup \rho_{0}=B_n\subset B(0,4)$, the first integral is bounded from above by a constant $a_f$. For the second integral we first note that 
\begin{equation*}
\log(1+\rho_{0}^{2}(z)) \leq  \log\left(1+\frac{4}{z_{0}^{2}}\right) \leq c\log\left( \frac{1}{z_0} \right), \qquad \text{where } z\in \mathbb{C},
\end{equation*}
for small $z_0$. For the rest of the integral we calculate 
\begin{equation*}
\int_{\mathbb{C}} \rho_{0}^2(z) dz \leq \sum_{j=1}^{n} \int_{B_j} \frac{4}{(2^{j-1}z_0)^2} dz \leq 4\pi n,
\end{equation*}
where $n$ was the smallest number for which $2^{n-1}z_0\geq 1$, and thus $n \leq c\log\left(\frac{1}{z_0}\right)$. Hence by combining the above estimates with \eqref{olkoon1} we obtain 
\begin{equation*}
M_{K_f}(\Gamma)\leq a_f+\frac{c}{p}\log^{2} \left( \frac{1}{z_0} \right),
\end{equation*}
which gives 
\begin{equation}\label{yksiviite}
M_{K_f}(\Gamma)\leq \frac{2c}{p}\log^{2} \left( \frac{1}{z_0} \right),
\end{equation}
for all $0<z_0<\epsilon_f$, for some $\epsilon_f>0$. \newline \newline
Next we will estimate $M(f(\Gamma))$ from below. Here we start with 
\begin{equation*}
M(f(\Gamma)) = \inf_{\rho \enskip \text{admissible}} \int_{\mathbb{C}} \rho^2(z) dz= \inf_{\rho \enskip \text{admissible}} \int_{0}^{2\pi} \int_{0}^{\infty} \rho^2(r,\theta)r dr d\theta,
\end{equation*}
and provide a lower bound for 
\begin{equation}\label{tulikinviite}
 \int_{0}^{\infty} \rho^2(r,\theta)r dr
\end{equation}
that holds for every direction $\theta\in [0,2\pi)$ and an arbitrary admissible $\rho$. We will estimate the integral \eqref{tulikinviite} from below by first finding $\left\lfloor\frac{\arg(f(z_0))-\arg(f(1))}{2\pi}\right \rfloor-1$ disjoint line segments $[x_i,y_i]\subset [0,1]$, for which one endpoint lies in $f(E)$ and the other in $f(F)$, and then using admissibility of $\rho$ to estimate its integral over these line segments. The main idea is to note that the paths $f(E)$ and $f(F)$ must cycle around the origin alternately, see figure \ref{FIG1} for illustration. \newline 
To see this, assume that the argument of the path $f(E)$ increases by $2\pi$ when $z$ moves from $t_0$ to $t_2$, where $z_0\leq t_2<t_0\leq 1$. As the mapping $f$ is a homeomorphism and the path $f(F)$ must contain the origin and points with big moduli, since $|f(t)|\to \infty$ when $t\to -\infty$, the path $f(F)$ must intersect the line segment $(f(t_2),f(t_0))$ at least once, let say at the point $f(t_1)$, where $-\infty<t_1<0$. We can choose $f(t_1)$ and $f(t_0)$ such that there are no points from the paths $f(E)$ and $f(F)$ in the line segment $(f(t_1),f(t_0))$. These line segments $[f(t_1),f(t_0)]$ are the ones we are looking for, and as the path $f(E)$ cycles around the origin $\left\lfloor\frac{\arg(f(z_0))-\arg(f(1))}{2\pi}\right \rfloor$ times we can define 
\begin{equation}\label{eräsviite}
n(z_0)=\left\lfloor\frac{\arg(f(z_0))-\arg(f(1))}{2\pi} \right \rfloor-1, 
\end{equation}
and find at least $n(z_0)$ disjoint line segments with desired endpoints. Thus we obtain, for every direction $\theta$, that
\begin{equation} \label{turha3}
 \int_{0}^{\infty} \rho^2(r,\theta)r dr \geq \sum_{i=1}^{n(z_0)} \int_{x_i}^{y_i} \rho^2(r,\theta)r dr,
\end{equation}
where 
\begin{equation*}\label{eitaas}
x_1<y_1<\cdots <x_{n(z_0)}<y_{n(z_0)}<1,
\end{equation*}
 and for every line segment $[x_i,y_i]$ the end points are in different sets $f(E),f(F)$. As every line segment $[x_i,y_i]$ belongs to the path family $f(\Gamma)$ we can estimate the integral over any line segment $[x_i,y_i]$, using the reverse H\"older inequality and admissibility of $\rho$, by
\begin{equation}\label{revhöl}
\int_{x_i}^{y_i} \rho^2(r,\theta)r dr \geq \left( \int_{x_i}^{y_i} \rho(r,\theta) dr \right)^2 \left( \int_{x_i}^{y_i} \frac{1}{r} dr \right)^{-1} \geq \frac{1}{\log\left( \frac{y_i}{x_i} \right)}.
\end{equation}
Combining this with \eqref{turha3} we obtain 
\begin{equation} \label{jatkatästä}
 \int_{0}^{\infty} \rho^2(r,\theta)r dr \geq \sum_{i=1}^{n(z_0)} \frac{1}{\log\left( \frac{y_i}{x_i} \right)}.
\end{equation}
To estimate this further we will use the following technical lemma.
\begin{lemma} \label{lemma}
Let $0<a<1$ and $n\in N$ be given, and let $a_i$ be positive numbers such that $a=a_0 <a_1<...<a_n=1$. Then 
\begin{equation}\label{summa1}
\sum_{i=1}^{n}\frac{1}{\log\left( \frac{a_i}{a_{i-1}} \right)} \geq \frac{n^2}{\log\left( \frac{1}{a} \right)}.
\end{equation}
\end{lemma}
\textit{Proof.} Choose the numbers $a_i$ such that we can find $a_{j-1}$, $a_j$ and $a_{j+1}$, for which $\frac{a_j}{a_{j-1}} \neq \frac{a_{j+1}}{a_j}$. By elementary calculations we see that replacing $a_j$ with $\bar{a}_j=\sqrt{a_{j-1}a_{j+1}}$ will decrease the value of the sum \eqref{summa1}. Thus, if we would know that the sum \eqref{summa1} attains its minimum we would be ready, as the numbers $a_i$ given by the condition $\frac{a_i}{a_{i-1}} = \frac{a_{i+1}}{a_i}$, for every $i$, would then have to give this minimum and can be calculated to satisfy the equality in \eqref{summa1}.\newline \newline
To show that the sum does attain its minimum we define 
\begin{equation*}
\Omega= \{x=(x_1,...,x_n)\in \mathbb{R}^n: x_i\geq 0 \enskip \text{for all $i$, and } x_1+\cdots +x_n=1-a\},
\end{equation*}
and note that it is compact. Then define the mapping $g$ on the set $\Omega$ by
\begin{equation*}
g(x)= \left\{
  \begin{array}{l l}
    \sum_{i=1}^{n} \frac{1}{\log\left( \frac{a+x_1+...+x_i}{a+x_1+...+x_{i-1}} \right)} & \quad \text{if $x_i>0$ for all $i$}\\
    \infty & \quad \text{if $x_i=0$ for some $i$ }
   
  \end{array} \right.
\end{equation*}
Note that if $x_i>0$ for every $i$ the sum \eqref{summa1} is coupled with $g(x)$ by the relation $x_i=a_i-a_{i-1}$. Let $x_0\in \Omega$ be such that $g(x_0)$ is finite. As $a$ is a fixed number there exists $\epsilon>0$ such that 
\begin{equation*}
\frac{1}{\log\left( \frac{a+x_1+...+x_i}{a+x_1+...+x_{i-1}} \right)}\geq g(x_0)
\end{equation*}
if $x_i< \epsilon$, for some $i$.  Hence $g(x)>g(x_0)$ for every $x\in \Omega$ such that $x_i<\epsilon$ for some $i$, and thus it is enough to look for the minimum in the set 
\begin{equation*}
\bar{\Omega}=\Omega\setminus \{x\in \mathbb{R}^n: x_i<\epsilon \enskip \text{for some $i$}\}.
\end{equation*}
As the set $\bar{\Omega}$ is compact and $g$ is continuous in the set $\bar{\Omega}$ the mapping $g$ attains its minimum in $\bar{\Omega}$ and hence also the sum \eqref{summa1} attains minimum, which finishes the proof.    \newline \newline
We return to proving our main theorem. We continue from \eqref{jatkatästä} by estimating 
\begin{equation}\label{kokoaminen}
\sum_{i=1}^{n(z_0)} \frac{1}{\log\left( \frac{y_i}{x_i} \right)} \geq \frac{1}{\log\left( \frac{y_1}{x_1}\right)}+ \sum_{i=2}^{n(z_0)-1} \frac{1}{\log\left( \frac{y_{i}}{y_{i-1}} \right)}+\frac{1}{\log\left( \frac{1}{y_{n(z_0)-1}}\right)} \geq \frac{n^2(z_0)}{\log\left( \frac{1}{x_1} \right)},
\end{equation}
where in the last inequality we have used Lemma \ref{lemma}. To estimate further from \eqref{kokoaminen} we must bound $x_1$ from below. To this end, use Theorem B from \cite{HK}, with the formulation of \eqref{modcont2}, on the modulus of continuity to get that 
\begin{equation}\label{viitteitälisää}
x_1\geq e^{-\frac{c}{p}\log^{2}\left( \frac{1}{z_0} \right)},
\end{equation}
when $z_0$ is small. Hence the estimate \eqref{kokoaminen} together with \eqref{viitteitälisää} gives
\begin{equation*}
\int_{0}^{\infty} \rho^2(r,\theta)r dr \geq \frac{p \thinspace n^2(z_0)}{c\log^2\left( \frac{1}{z_0} \right)},
\end{equation*} 
for an arbitrary direction $\theta$ and small $z_0$. Thus we have for these points $z_0$ that 
\begin{equation}\label{lasku1}
M(f(\Gamma)) = \inf_{\rho \enskip \text{admissible}} \int_{0}^{2\pi} \int_{0}^{\infty} \rho^2(r,\theta)r dr d\theta \geq \frac{c p \thinspace n^2(z_0)}{\log^2 \left( \frac{1}{z_0} \right)}.
\end{equation}
We now have the estimates \eqref{yksiviite} and \eqref{lasku1} for the moduli $M_{K_f}(\Gamma)$ and $M(f(\Gamma))$, when the point $z_0$ is small. Then we use the inequality \eqref{A} with these estimates to obtain 
\begin{equation*}
\frac{cp \thinspace n^2(z_0)}{\log^2 \left( \frac{1}{z_0} \right)} \leq \frac{c}{p}\log^2 \left( \frac{1}{z_0} \right),
\end{equation*}
and see
\begin{equation*}
n(z_0) \leq \frac{c}{p} \log^2 \left( \frac{1}{z_0} \right).
\end{equation*}
Thus, due to the equality \eqref{eräsviite}, we have that
\begin{equation*}
\arg(f(z_0))-\arg(f(1)) \leq  \frac{c}{p} \log^2 \left( \frac{1}{z_0} \right),
\end{equation*}
when $z_0<\epsilon_f$. And since $\arg(f(1))$ is finite for any fixed branch of the argument this proves Theorem \ref{Theorem2}, and hence also its renormalization Theorem \ref{Theorem1}. \newline \newline
We will next show that Theorem \ref{Theorem1} is sharp up to the constant $c$ by proving that, given arbitrary parameters $c_1$ and $c_2$ satisfying the condition \eqref{vakioyht}, the mapping 
\begin{equation}\label{esimerkki}
h(z)=  \left\{
  \begin{array}{l l}
    \frac{z}{|z|}|z|^{c_1\log\left( \frac{1}{|z|} \right)-ic_2\log\left( \frac{1}{|z|} \right)} & \quad \text{if $z\in B\left(0,\frac{1}{2}\right)$}\\
     \frac{z}{|z|}|z|^{c_1\log(2)-ic_2 \log(2)} & \quad \text{if } z\in \mathbb{D}\setminus B\left(0,\frac{1}{2}\right)\\
    z & \quad \text{if $z\notin \mathbb{D}$ }
  \end{array} \right.
\end{equation}
has $p$-exponentially integrable distortion and that we can choose the parameters $c_1>0$ and $c_2\in \mathbb{R}$ such that the rotation is close to the bound given by Theorem \ref{Theorem1}. First note, that it is enough to calculate the distortion function $K_h(z)$ in the ball $ B\left(0,\frac{1}{2}\right)$, as otherwise it is bounded by some constant. To do this compute the complex derivatives 
\begin{equation}\label{zder}
h_z(z)= \frac{e^{-c_1\log^{2}\left( \frac{1}{|z|} \right)+ic_2\log^{2}\left( \frac{1}{|z|} \right)}}{|z|} \left( \frac{1}{2}+c_1 \log\left( \frac{1}{|z|}\right)-ic_2\log\left( \frac{1}{|z|} \right)  \right),
\end{equation}
and
\begin{equation}\label{barzder}
h_{\bar{z}}(z)= \frac{\sqrt{z} e^{-c_1\log^{2}\left( \frac{1}{|z|} \right)+ic_2\log^{2}\left( \frac{1}{|z|} \right)}}{\bar{z}^{\frac{3}{2}}} \left( -\frac{1}{2}+c_1 \log\left( \frac{1}{|z|} \right)-ic_2\log\left( \frac{1}{|z|} \right) \right).
\end{equation}
Then calculate the modulus of the Beltrami coefficient
\begin{equation}\label{mu}
|\mu_h(z)|=\frac{|h_{\bar{z}}(z)|}{|h_z(z)|}=\frac{\sqrt{\left( c_1 \log\left( \frac{1}{|z|} \right) -\frac{1}{2}\right)^2+c_{2}^{2}\log^{2} \left( \frac{1}{|z|} \right)}}{\sqrt{\left( c_1 \log\left( \frac{1}{|z|} \right) +\frac{1}{2}\right)^2+c_{2}^{2}\log^{2} \left( \frac{1}{|z|} \right)}}.
\end{equation}
And finally use the equation $K_h(z)=\frac{1+|\mu_h(z)|}{1-|\mu_h(z)|}$ with \eqref{mu} to obtain
\small \begin{equation*}
\begin{split}
K_h(z) & = \frac{\sqrt{\left( c_1 \log\left( \frac{1}{|z|} \right) +\frac{1}{2}\right)^2+c_{2}^{2}\log^{2} \left( \frac{1}{|z|} \right)}+\sqrt{\left( c_1 \log\left( \frac{1}{|z|} \right) -\frac{1}{2}\right)^2+c_{2}^{2}\log^{2} \left( \frac{1}{|z|} \right)}}{\sqrt{\left( c_1 \log\left( \frac{1}{|z|} \right) +\frac{1}{2}\right)^2+c_{2}^{2}\log^{2} \left( \frac{1}{|z|} \right)}-\sqrt{\left( c_1 \log\left( \frac{1}{|z|} \right) -\frac{1}{2}\right)^2+c_{2}^{2}\log^{2} \left( \frac{1}{|z|} \right)}}\\
& = \frac{\xi(|z|) 4(c_{1}^{2}+c_{2}^{2})\log^{2}\left( \frac{1}{|z|} \right)}{2c_1 \log\left( \frac{1}{|z|} \right)}=\xi(|z|) \left( 2c_1+\frac{2c_{2}^{2}}{c_1} \right) \log\left( \frac{1}{|z|} \right),
\end{split}
\end{equation*} \normalsize
where $\xi$ is bounded, for $z\in  B\left(0,\frac{1}{2}\right)$, and $\xi(|z|)\to 1$ as $|z|\to 0$. Thus $K_h(z)$ is $p$-exponentially integrable when 
\begin{equation*}
2c_1+\frac{2c_{2}^{2}}{c_1} < \frac{2}{p}.
\end{equation*}
From this we obtain that  
\begin{equation}\label{lklk}
0<c_1<\frac{1}{p},
\end{equation}
and
\begin{equation*}
|c_2| < \sqrt{c_1\left( \frac{1}{p}-c_1 \right)} \leq \frac{1}{2p},
\end{equation*}
where in the last inequality we use the bound \eqref{lklk} for $c_1$. This shows that given an arbitrary $\epsilon>0$ we can choose the parameters $c_1$ and $c_2$ such that the mapping \eqref{esimerkki} satisfies
\begin{equation*}
|\Im \log (h(r))|= \frac{1-\epsilon}{2p} \log^{2}\left( \frac{1}{r} \right)
\end{equation*}
for all $0<r<\frac{1}{2}$. Clearly the mappings $h$ are sense-preserving homeomorphisms. From the calculations \eqref{zder} and \eqref{barzder}, together with the definition \eqref{esimerkki}, we see that they lie in the Sobolev space $W^{1,1}_{\text{loc}}(\mathbb{C})$, and thus also have locally integrable Jacobian. Hence these mappings have finite distortion and thus prove that Theorem \ref{Theorem1} is sharp up to the constant $c$, and additionally that the rotation 
\begin{equation*}
\frac{1-\epsilon}{2p}
\end{equation*}
can be attained. \newline \newline
Regarding the remarks on the relation between stretching and rotation we first calculate the distortion of the mappings \eqref{esimkiertoeiveny}, defined by 
\begin{equation*}
h_0(z)=  \left\{
  \begin{array}{l l}
    ze^{-ic_2\log^{\frac{3}{2}}\left( \frac{1}{|z|} \right)} & \quad \text{if $z\in \mathbb{D}$}\\
    
    z & \quad \text{if $z\notin \mathbb{D}$ }
  \end{array} \right. 
\end{equation*}
in a similar manner as for the mappings \eqref{esimerkki}, and obtain that 
\begin{equation*}
K_{h_0}(z)=1+\xi(|z|)\frac{9}{4}c_{2}^{2}\log\left( \frac{1}{|z|} \right),
\end{equation*}
where $\xi$ is bounded and $\xi(|z|)\to 1$ as $|z|\to 0$. This proves that $h_0$ has $p$-exponentially integrable distortion when 
\begin{equation*}
\frac{9c_{2}^{2}}{8}< \frac{1}{p},
\end{equation*}
and with a similar reasoning as for the mappings \eqref{esimerkki} we see that $h_0$ is a mapping of finite distortion. Moreover, since $|h_0(z)|=|z|$, we see that $\arg(h_0(r))$ grows like $\log^{\frac{3}{2}}\left( \frac{1}{|h_0(r)|} \right)$ as $r\to 0$. \newline \newline 
However, the pointwise stretching of $f$ does bound the pointwise rotation even in the case of mappings with exponentially integrable distortion, but the relation between stretching and rotation is different from the quasiconformal case. This follows from \eqref{viitteitälisää}, where we use the modulus of continuity. If in \eqref{viitteitälisää} we would instead of 
\begin{equation*}
x_1\geq e^{-\frac{c}{p}\log^{2}\left( \frac{1}{z_0} \right)} 
\end{equation*}
assume 
\begin{equation}\label{jj}
x_1\geq e^{-\frac{c}{p}\log^{1+\alpha}\left( \frac{1}{z_0} \right)},
\end{equation}
where $\alpha\in (-1,1)$, and continue after \eqref{viitteitälisää} as in the proof of Theorem \ref{Theorem2} we would obtain
\begin{equation}\label{hhhhh}
\arg(f(z_0))-\arg(f(1)) \leq \frac{c_s}{p}\log^{\frac{3+\alpha}{2}}\left( \frac{1}{z_0} \right),
\end{equation}
where $c_s$ is a constant depending on the constant $c$ chosen in \eqref{jj}. This shows how stretching bounds rotation from above. To see that this relation between stretching and rotation is optimal, again up to the constant $c_s$, we present the mappings
\begin{equation}
h_{\alpha}(z)=  \left\{
  \begin{array}{l l}
    \frac{z}{|z|}|z|^{c_1\log^{\alpha}\left( \frac{1}{|z|} \right)-ic_2\log^{\frac{1+\alpha}{2}}\left( \frac{1}{|z|} \right)} & \quad \text{if $z\in B\left(0,\frac{1}{2}\right)$}\\
    \frac{z}{|z|}|z|^{c_1\log^{\alpha}(2)-ic_2 \log^{\frac{1+\alpha}{2}}(2)} & \quad \text{if } z\in  \mathbb{D}\setminus B\left(0,\frac{1}{2}\right) \\
    z & \quad \text{if $z\notin \mathbb{D}$ }
  \end{array} \right. 
\end{equation}
where $c_1>0$ and $c_2\in \mathbb{R}$. These mappings can be checked to be mappings of finite distortion with $p$-exponentially integrable distortion when
\begin{equation}\label{jjjj}
\frac{c_{2}^{2}\left(\frac{3+\alpha}{2} \right)^2}{2c_1(1+\alpha)}< \frac{1}{p},
\end{equation}
with a similar calculation as for the mappings \eqref{esimerkki}. From \eqref{jjjj} we see that $c_2$ can be arbitrary big if we choose sufficiently big $c_1$, and hence the constant $c_s$ truly does depend on the stretching and there is no absolute constant $c$ for which \eqref{hhhhh} would hold.

\end{document}